\documentclass[11pt]{article}
\usepackage[a4paper]{geometry}

\usepackage{amsmath}
\usepackage{amsfonts}
\usepackage{graphicx}
\usepackage{hyperref}
\DeclareMathAlphabet{\altmathcal}{OMS}{cmsy}{m}{n}

\title{Time-reversible dynamics in a system of two coupled active rotators}
\author{
  Oleksandr Burylko$^{1,2,3}$, Matthias Wolfrum$^{3}$, \\
  Serhiy Yanchuk$^{1,4}$, and Jürgen Kurths$^{1,5}$\\[0.7cm]
  $^{1}$Potsdam Institute for Climate Impact Research, 14437 Potsdam, Germany\\
  $^{2}$Institute of Mathematics, National Academy of Sciences of Ukraine,\\ Tereshchenkivska Str. 3, 01024 Kyiv, Ukraine\\
  $^{3}$Weierstrass Institute for Applied Analysis and Stochastics, \\
  Mohrenstrasse 39, 10117 Berlin, Germany\\
  $^{4}$Institute of Mathematics, Humboldt University Berlin, 12489 Berlin, Germany\\
  $^{5}$Department of Physics, Humboldt University Berlin, 12489 Berlin, Germany
}

\begin{document}





\maketitle

\begin{abstract}
We study two coupled active rotators with Kuramoto-type coupling and focus our attention to specific transitional regimes where the coupling is neither attractive nor repulsive. We show that certain such  situations at the edge of synchronization can be characterized by the existence of a time-reversal symmetry of the system. We identify two different cases with such a time-reversal symmetry. The first case is characterized by a non-reciprocal attractive/repulsive coupling. The second case is a reciprocal coupling exactly at the edge between attraction and repulsion. We give a detailed description of possible different types of dynamics and bifurcations for both cases. In particular, we show how the time-reversible coupling can induce both oscillation death and oscillation birth to the active rotators. Moreover, we analyse the coexistence of conservative and dissipative regions in phase space, which is a typical feature of systems with a time-reversal symmetry. We show also, how perturbations breaking the time-reversal symmetry and destroying the conservative regions can lead to complicated types of dissipative dynamics such as the emergence of long-period cycles showing a bursting-like behavior.
\end{abstract}


\section{\label{sec: Introduction}Introduction}

Collective dynamics of weakly interacting oscillatory systems can be effectively described by coupled phase oscillators \cite{Kuramoto1984,Winfree2001,Hoppensteadt1997}. 
The classical Kuramoto system of coupled phase oscillators is based on the assumption that without coupling each subsystem has oscillatory (periodic) dynamics and hence, for weak coupling,  can be reduced to the simple phase equation $\dot \phi_j = \omega_j$ with some internal frequency $\omega_j$. 
In particular, it has been used extensively for the study of various forms of synchronization  \cite{Acebron2005,RodriguesPeronJiEtAl2016,Pikovsky2001}. 
Systems of coupled active rotators have already been introduced by Shinomoto and Kuramoto in 1986 \cite{Shinomoto1986} to study a more general class of interacting units, where each unit is governed by a non-homogeneous oscillator $\dot \phi_j= \omega_j - a_j \cos \phi_j$. 
In particular such units undergo for $|a_j|=|\omega_j|$  a so-called SNIC (saddle-node on invariant circle) bifurcation such that the oscillator is transformed into an excitable unit. In this sense, coupled active rotators provide a substantial extension compared to the classical Kuramoto system of phase oscillators and  are suitable for the modeling of collective dynamics of neuronal and, in general, excitable systems. The active rotator of this form is also known as the theta-neuron model \cite{Bick2020,Laing2020,So2014,Ashwin2016a,Zheng2018}, and it is equivalent to the quadratic integrate-and-fire neuron \cite{Ermentrout1986, Montbrio2015}.
Systems of coupled active rotators and their extensions have been also studied in~\cite{Sakaguchi_1988, Kurrer_1995, Park_1996, Tessone_2007, Zaks_2003, Lafuerza_2010, Ionita_2014, OKeeffe_2016, Zaks2016, Dolmatova_2017, Klinshov_2019, 
Franovic_2021, Klinshov_2021, Wolfrum_2021, Zaks_2021, Bacic2018, Franovic2020, Yanchuk2019}. 

The onset of various forms of synchronization and collective dynamics is usually studied in the context of attractive global or non-local coupling. However, many interesting and unexpected dynamical effects can be observed close to the transition from attractive to repulsive coupling \cite{Kuramoto-Battogtokh} and also  units with different types of non-reciprocal coupling can lead to new dynamical phenomena \cite{Hong-Strogatz-PRL,Hong-Strogatz-PRE2011}.
In this work, we consider a minimal network motif of  two coupled active rotators with Kuramoto-type coupling and focus our attention to  specific transitional regimes where the coupling is neither attractive nor repulsive.
 It turns out that certain situations at the edge of synchronization can be characterized by an additional structural property of the system, the existence  of a time-reversal symmetry. Systems with this property are known to exhibit rich and unexpected dynamical behavior and have been studied extensively from a mathematical point of view \cite{Sevryuk1986, Topaj2002, Burylko2018, Politi1986, Lamb1998lin,ASH16a}. A specific feature of such systems is the possibility of a coexistence of  regions with conservative dynamics (e.g. families of neutrally stable closed orbits) with dissipative regions in phase space. 
 
 In our setting of two coupled rotators we identify two different cases with such a time-reversal symmetry. The first case is characterized by a non-reciprocal coupling where one oscillator couples attractive and the other repulsive. The second case is a reciprocal coupling with a Kuramoto type coupling exactly at the edge between attractive and repulsive  coupling. 
For both cases we describe in detail the different dynamical scenarios and the  bifurcation transitions between them. In particular, we show how the coupling can lead to coexistence of rotations in opposite directions, to the birth and death of oscillations, and to the coexistence of a dissipatively stable synchronous equilibrium  with regions of conservative oscillatory motions in the form of  both rotations and  librations.
In section \ref{sec:pert}, we also consider the influence of generic perturbations. Such perturbations can induce a drift along the families of periodic solutions in the conservative regions of the reversible regime. We show that in some cases  this can lead to  long-period limit cycles  with a bursting-like dynamics. 
Additionally, we study the effects of the higher Fourier modes  that can  lead to even higher multistablity of conservative and dissipative regions. 


%
%

A general system of two coupled rotators has the form
\begin{eqnarray}
\dot \phi_{1} = f_{1}(\phi_{1})+g_{1}(\phi_{1}-\phi_{2}),
\label{eq: General_1}\\
\dot \phi_{2} = f_{2}(\phi_{2})+g_{2}(\phi_{2}-\phi_{1}),
\label{eq: General_2}
\end{eqnarray}
where $\phi_{1}$, $\phi_{2}\in\mathbb{T}^{1}=\mathbb{R}/2\pi \mathbb{Z}$ are phase variables, and the local dynamics $f_{1,2}$ as well as the coupling functions $g_{1,2}$ are smooth and $2\pi$-periodic. 
We mainly restrict ourselves to the case
\begin{eqnarray}
\dot \phi_{1} = \omega_{1} + a_{1}\cos\phi_{1} + \kappa_{1}\sin(\phi_2 - \phi_1 + \alpha),
\label{eq:rot1} \\
\dot \phi_{2} = \omega_{2} + a_{2}\cos\phi_{2} + \kappa_{2}\sin(\phi_1 -\phi_2 + \alpha),
\label{eq:rot2}
\end{eqnarray}
where both the local dynamics and the coupling functions contain only the leading Fourier component. In this way we get the  natural frequencies  $\omega_{i}$,  the phase inhomogeneities $a_{i}$ the coupling strengths $\kappa_{i}$ and  the phase shift $\alpha$ as parameters. More complicated functions $f_i$ and $g_i$ will be also shortly discussed, and they will be specified at the corresponding places. 

The inhomogeneity $a_i$ is an important ingredient of the system, since otherwise the dynamics is very simple.  Indeed, if $a_{1}=a_{2}=0$, we obtain two coupled oscillators of Kuramoto-Sakaguchi type \cite{Strogatz1994}, which have a
phase-shift symmetry $(\phi_{1},\phi_{2})\mapsto(\phi_{1}+\delta,\phi_{2}+\delta)$ for any $\delta\in\mathbb{T}^{1}.$ 
As a result, the system can be reduced to a single equation for the phase difference $\psi=\phi_{1}-\phi_{2}$ 
\[
\dot \psi
=\Delta-A\cos(\psi-\sigma),
\]
where $\Delta=\omega_{1}-\omega_{2}$, $\tan\sigma=\frac{\kappa_{1}+\kappa_{2}}{\kappa_{2}-\kappa_{1}}\cot\alpha$, and
$$A=\sqrt{\left(\kappa_{2}-\kappa_{1}\right)^{2}\sin^{2}\sigma+\left(\kappa_{1}+\kappa_{2}\right)^{2}\cos^{2}\sigma}.$$
This is again an active rotator with stable and unstable equilibria for $|\Delta/A|<1$ and the SNIC bifurcation for $\Delta/A=\pm1$. 
The stable and unstable equilibria for the system in the phase differences correspond to stable and unstable phase-locked limit cycles for the original two-dimensional Kuramoto-Sakaguchi system. 
For $|\Delta/A|>1$, the phase-locking is lost, and the system of two coupled Kuramoto-Sakaguchi oscillators possesses families of neutral periodic or quasi-periodic orbits depending on the relationship between $\omega_{1}$ and $\omega_{2}$.

The dynamics of system (\ref{eq:rot1})--(\ref{eq:rot2}) becomes more complicated when the inhomogeneity $a_{i}$ is present. The phase-shift symmetry is broken and the transition to the excitable regime of the single unit induces new dynamical regimes. 
As we will see, the dynamics are particularly rich are in the cases of time-reversible coupling as discussed in the following section. 

%
%
\section{Time-reversible dynamics of two coupled oscillators \label{sec:time-reversible}}

\subsection{What is time-reversibility?}
A system $\dot x = F(x)$ has a time-reversal symmetry \cite{Sevryuk1986, Topaj2002, Burylko2018, Politi1986, Lamb1998} if there exists an involution $\altmathcal{R}$ of the phase space $X$ satisfying
\begin{equation}
F(\altmathcal{R}(x))=-\altmathcal{R}(F(x))\label{eq: Time-rev_symm}
\end{equation}
and $\altmathcal{R}^{2}=\mathrm{Id}$, with $\mathrm{Id}$ being the identity transformation.
The existence of such a time-reversing symmetry action $\altmathcal{R}$, which is typically linear or affine,  implies that for a solution $x(t)$ also  $\altmathcal{R}x(-t)$ is a solution. An important role for the characterization of the dynamics of a reversible systems plays the subspace 
\begin{equation}
\mathrm{Fix}\altmathcal{R}=\left\{x\in X: \altmathcal{R}(x)=x\right\}.
\label{eq: FixR}
\end{equation}
In contrast to invariant subspaces of symmetries without time reversal, this subspace is not dynamically invariant. Instead, a trajectory can  cross $\mathrm{Fix}\altmathcal{R}$, which then implies that the whole trajectory is mapped by $\altmathcal{R}$ onto a time reversed copy of itself. In this way one can  distinguish between intersecting trajectories connecting an attractor-repellor pair related by  $\altmathcal{R}$, and trajectories intersecting more than once, which induces locally conservative dynamics. This coexistence of conservative and dissipative dynamics in different regions of the phase space is a typical property of systems with time-reversibility \cite{Politi1986,Golubitsky-Krupa-Lim,Topaj-Pikovsky,BMWY_SIAM_2018} that distinguishes them from generic dissipative dynamical systems. 

In systems with time-reversal symmetry one has to distinguish between equilibria within and outside  $\mathrm{Fix}\altmathcal{R}$. In the first case, an equilibrium has to have the same number of stable and unstable directions, since these are related by $\altmathcal{R}$. In the second case equilibria come in pairs, related by $\altmathcal{R}$, with opposite stability properties. Moreover, there can be bifurcations with a spontaneous symmetry breaking, where from a branch of equilibria within  $\mathrm{Fix}\altmathcal{R}$ a branch containing pairs of equilibria  outside $\mathrm{Fix}\altmathcal{R}$ bifurcates, see e.g. \cite{Lim-McComb_1998,Broer_1998,Lerman-Turaev2012}. Due to the possibility of locally conservative dynamics, reversible systems can have structurally stable homoclinic orbits and heteroclinic cycles, which  together with their specific bifurcations have been studied extensively, see \cite{Knobloch1997,Wagenknecht_2005,Homburg_2006}.  

\subsection{Reversible cases in the system of coupled rotators\label{subsec: Revers General}}
We identify two cases, for which system (\ref{eq: General_1})--(\ref{eq: General_2}) is time-reversible. 
First, note that the single rotator has a time reversal symmetry 
$$(\phi,\,t)\longmapsto(-\phi_,\,-t)$$
as soon as $f$ is an even function. A coupled system of two identical such units, i.e. with
$$f_{1}(\phi)=f_{2}(\phi) =f(\phi) ,\quad  f(\phi)=f(-\phi),$$ can become time-reversible in two different ways. 
Case (I) is characterized by 
an anti-reciprocal coupling with an odd coupling function
\begin{equation}
g_{1}(\phi)=-g_{2}(\phi)=g(\phi),\quad g(-\phi)=-g(\phi).\label{eq: Cond_g}
\end{equation}


The second time-reversible case (II) appears for if the coupling functions are identical and even.
\begin{equation} g_{1}(\phi)=g_{2}(\phi)=g(\phi),\quad  g(\phi)=g(-\phi),
\end{equation}
which corresponds to a conservative coupling at the edge between attraction and repulsion. 
In both cases,  the time-reversible symmetry is given by the action
\begin{equation}
\altmathcal{R}:\ (\phi_{1},\,\phi_{2},\,t)\longmapsto(-\phi_{2},\,-\phi_{1},\,-t)
\label{eq: R_Sym}
\end{equation}
with the subspace
\begin{equation}
\mathrm{Fix}\altmathcal{R}=\left\{ (\phi_{1},\,\phi_{2}):\,\phi_{1}=-\phi_{2}\right\}.
\label{eq: FixR10}
\end{equation}
For the system (\ref{eq:rot1})--(\ref{eq:rot2}) of active rotators with Kuramoto-Sakaguchi type coupling we obtain for case (I) with anti-reciprocal and odd coupling the system
\begin{eqnarray}
\dot \phi_{1} = \omega + a\cos\phi_{1}-\kappa\sin(\phi_{1}-\phi_{2}),\label{eq:CaseI_1} \\
\dot \phi_{2} = \omega + a\cos\phi_{2}+\kappa\sin(\phi_{2}-\phi_{1}), \label{eq:CaseI_2}
\end{eqnarray}
while in case (II) with even and reciprocal coupling we get 
\begin{eqnarray}
\dot \phi_{1} = \omega + a\cos\phi_{1}-\kappa\cos(\phi_{1}-\phi_{2}),\label{eq:CaseII_1}\\
\dot \phi_{2} = \omega + a\cos\phi_{2}-\kappa\cos(\phi_{2}-\phi_{1}).\label{eq:CaseII_2}
\end{eqnarray}
In the following Sec.~\ref{sec:dynamics-reversible}, we describe the dynamics and bifurcations in the above two cases. 

%
%

\section{Time-reversible dynamics of the coupled rotator model \label{sec:dynamics-reversible}}
\subsection{Case (I): coupled rotators with anti-reciprocal coupling}
We first consider the case (I) reversible system (\ref{eq:CaseI_1})--(\ref{eq:CaseI_2}).
This system possesses additional symmetries that involve  parameters; these symmetries are generated by the actions
\begin{align}
   \gamma_{1}: & \quad (\phi_{1},\,\phi_{2},\omega,\,\,t)\,\longmapsto(\phi_{2}+\pi,\,\phi_{1}+\pi,\,-\omega,\,-t),\\ 
   \gamma_{2}: & \quad  (\phi_{1},\,\phi_{2},\kappa,\,\,t) \;\longmapsto(\phi_{2}+\pi,\,\phi_{1}+\pi,\,-\kappa,\,t),\\
   \gamma_{3}: & \quad  (\phi_{1},\,\phi_{2},a,\,\,t) \,\; \longmapsto(\phi_{1}+\pi,\,\phi_{2}+\pi,\,-a,\,t).
\end{align}

Note that $\gamma_1$ induces for $\omega=0$ a second time-reversing symmetry action, while $\gamma_{2,3}$ for $\kappa=0$ and $a=0$, respectively, induce $\mathbb{Z}_2$-symmetries without time reversal. 
As a result of the parametric symmetries $\gamma_{1,2,3}$ the resulting  bifurcation diagrams will be mirror symmetric with respect all the parameters $\omega,\,\kappa,\,a$.
Also the synchrony subspace $\phi_1=\phi_2$ is flow invariant for system (\ref{eq:CaseI_1})--(\ref{eq:CaseI_2}). However, this invariance is not induced by a symmetry of the system, but the diffusive nature of the coupling. 
\begin{figure}
\centering{}\includegraphics[width=\textwidth]{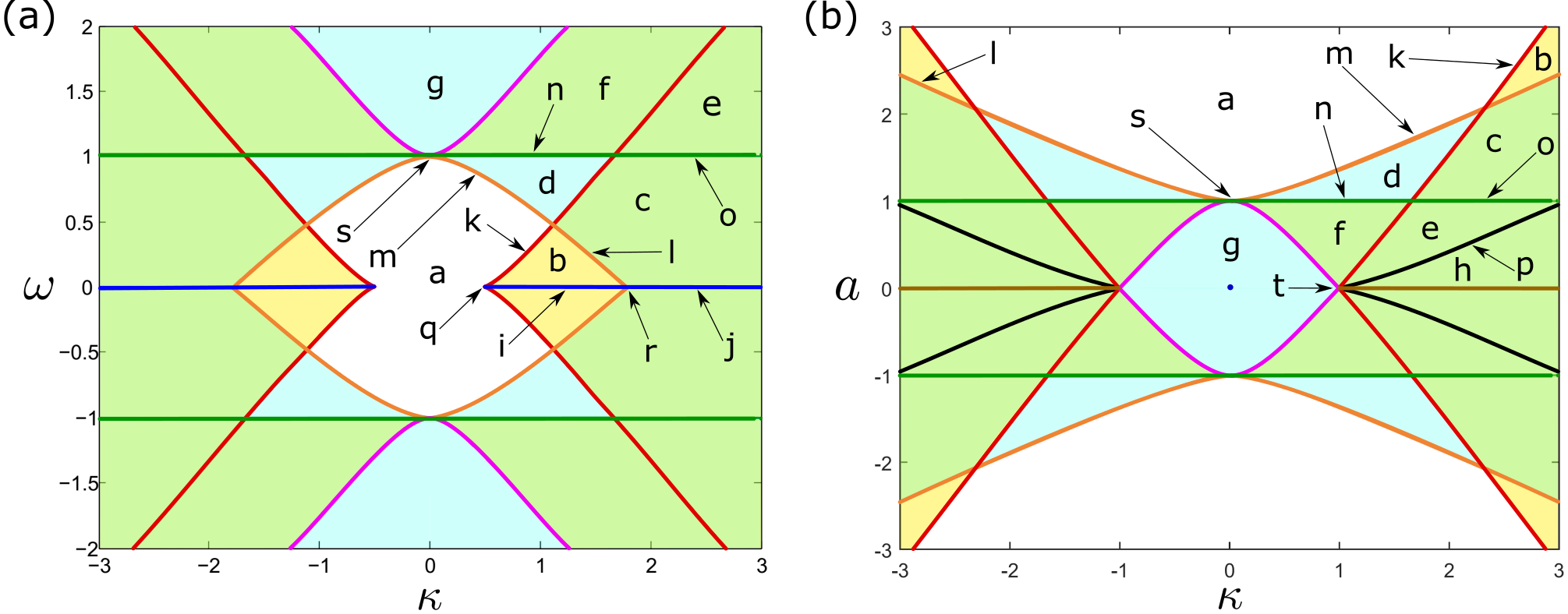}
\caption{\label{Fig: Bif_Diagr_S2} 
 Bifurcation diagrams for system (\ref{eq:CaseI_1})--(\ref{eq:CaseI_2}).
(a) -- parameter plane $(\kappa,\omega)$  with fixed $a=1$. 
(b) -- parameter plane  $(\kappa,a)$ with fixed $\omega=1$.
Bifurcation curves: red/magenta -- saddle-center  bifurcation, 
green -- reversible pitchfork bifurcation,
orange/black -- heteroclinic saddle-saddle connections.
Global bifurcations induced by a second time-reversal symmetry (blue) and by the phase shift invariance (brown). Structurally stable phase portraits from the different parameter regions with  labels (a)--(h) are given in Fig.~\ref{Fig:Phase_Portr_1}, phase portraits at the bifurcation curves with labels (i)--(p) are given in Fig.~\ref{Fig:Phase_Portr_2} and phase portraits of the codimension-two points with labels (q)--(t) in Fig.~\ref{Fig:Phase_Portr_3}. The colors of the regions indicate purely dissipative  dynamics (white), 
mixed type with a dissipative and a libration region (yellow), different types without libration regions (blue), different types  with librations and rotations (green).
}
\end{figure}

The regions in the bifurcation diagrams in Fig.~\ref{Fig: Bif_Diagr_S2} correspond to qualitatively different structurally stable phase portraits. Panel (a) shows the parameter plane  $(\kappa,\omega)$ with fixed $a=1$ and panel (b)  the plane $(\kappa,a)$  with fixed $\omega=1$. 
Note that one can fix $a=1$  without loss of generality as soon as $a\ne 0$. 
To study the situation  in a vicinity of $a=0$, we can fix $\omega=1$, instead. 
The region $a\approx 0$ is interesting from the point of view of perturbing the Kuramoto system 
with its phase shift symmetry to a rotator system with inhomogeneous rotation speed.  
Note that the  diagram in panel (b) can be obtained from the diagram in panel (a) by  the transformation
\begin{equation}
\label{eq:transform}
(\kappa,\omega) \mapsto (\kappa/\omega,1/\omega).
\end{equation}
This explains, why the blue and brown bifurcation curves, which are mapped from or to infinity, are present only in one of the two diagrams. 
The black  line in panel (b) is not visible in panel (a) because its preimage lies outside the plotted region.
In the remaining part of this section we will describe in detail the different types of bifurcations  indicated in the two diagrams by curves of different color.

\begin{figure*}
\centering{}
\includegraphics[width=1\textwidth]{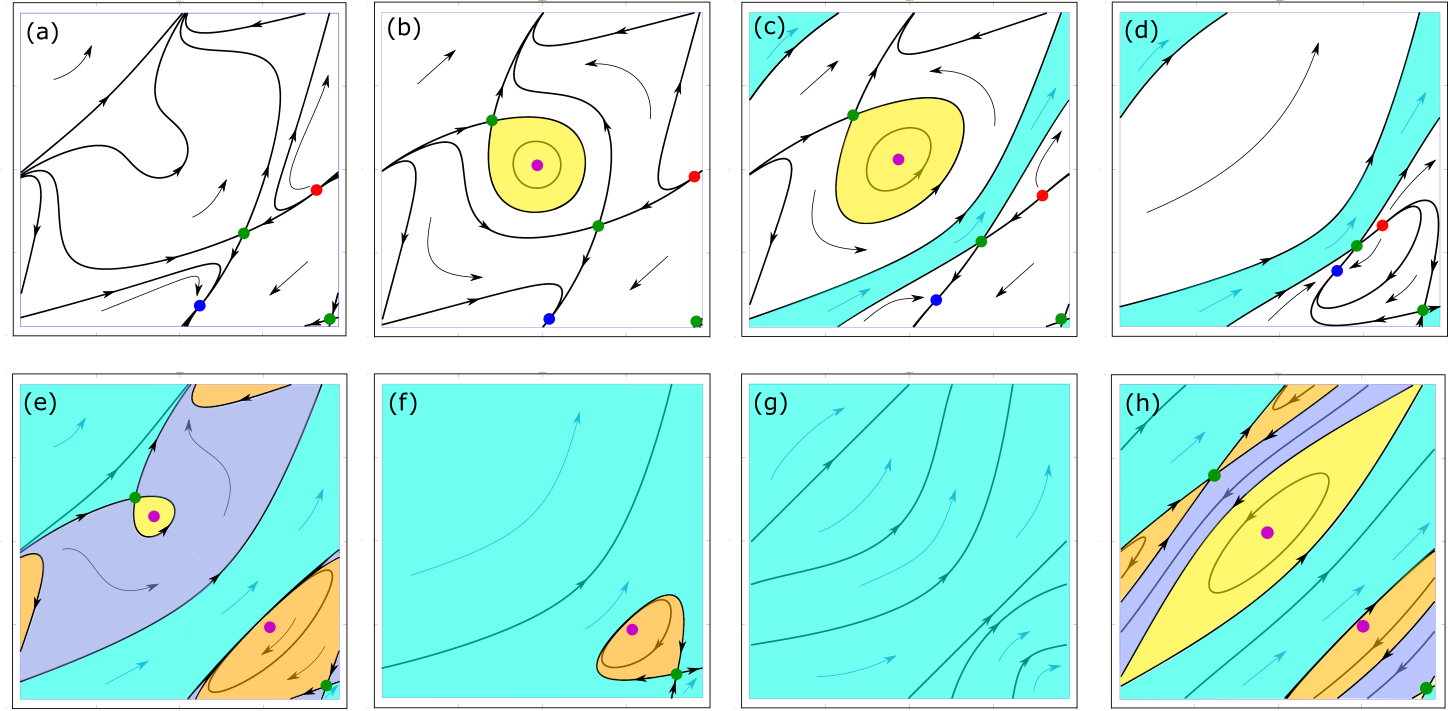}
\caption{
Different types of structurally stable phase portraits of system (\ref{eq:CaseI_1})--(\ref{eq:CaseI_2}) in case (I) of time reversibility  with anti-reciprocal coupling. 
The squared regions show the torus $(\phi_{1},\phi_{2})\in[-\pi/2,3\pi/2]\times[-3\pi/2,\pi/2]$.
Parameters in panels (a)--(h) are chosen from the correspondingly marked regions from the bifurcation diagrams in Fig.~\ref{Fig: Bif_Diagr_S2}. Conservative regions are colored as follows: Orange -- librations with clockwise motion, yellow -- librations with anti-clockwise motion, blue/cyan -- rotations. Dissipative regions are white.
Fixed points are colored as follows: red --- source, blue --- sink, green --- saddle, magenta --- center.
}\label{Fig:Phase_Portr_1} 
\end{figure*}
\begin{figure*}
\centering{}
\includegraphics[width=1\textwidth]{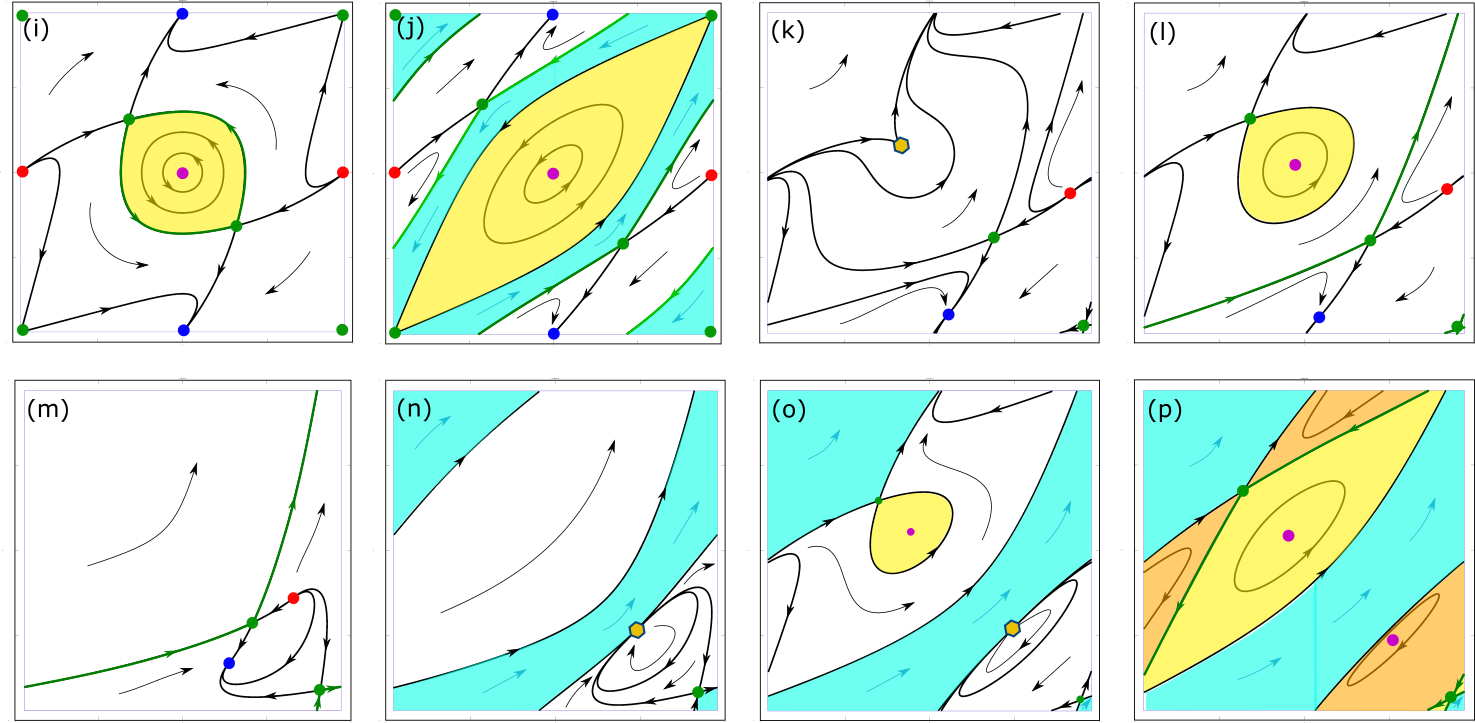}
\caption{\label{Fig:Phase_Portr_2} 
Examples of structurally unstable phase portraits of system (\ref{eq:CaseI_1})--(\ref{eq:CaseI_2}). 
Parameters are chosen from the correspondingly marked bifurcation curves in Fig.~\ref{Fig: Bif_Diagr_S2}. Degenerate equilibria (orange hexagons): (k) -- saddle-center bifurcation, (n), (o) -- reversible pitchfork. Structurally unstable saddle-saddle connections are marked in green in panels (i), (j), (l), (m), and (p). Colors of invariant  regions and fixed points as in Fig.~\ref{Fig:Phase_Portr_1}.
}
\end{figure*}
\begin{figure*}
\centering{}
\includegraphics[width=1\textwidth]{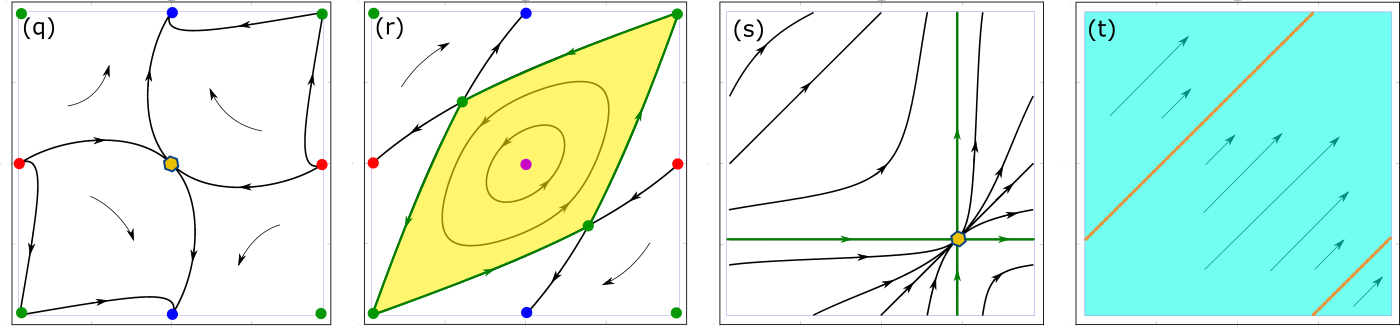}
\caption{\label{Fig:Phase_Portr_3} 
Phase portraits at the  correspondingly marked codimension-two points of system (\ref{eq:CaseI_1})--(\ref{eq:CaseI_2}), shown in Fig.~\ref{Fig: Bif_Diagr_S2}. Degenerate equilibria are marked with hexagons (orange). Structurally unstable saddle-saddle connections are marked in green. Colors of invariant  regions and fixed points as in Fig.~\ref{Fig:Phase_Portr_1}.}
\end{figure*}

Examples of the different generic phase portraits are depicted in Fig.~\ref{Fig:Phase_Portr_1}.  In Fig.~\ref{Fig:Phase_Portr_2} we show examples of structurally unstable phase portraits on the different bifurcation curves and  Fig.~\ref{Fig:Phase_Portr_3} gives the phase portraits at the codimension-two points. The different dynamical regimes are distinguished by the number and type of the fixed points and also by homoclinic and heteroclinic connections that organize the dissipative and conservative regions. 

The fixed points have the following general properties:
\begin{itemize}
\item Depending on the  parameter values, the system has up to six fixed points.
\item There can be up to four fixed points in  $\mathrm{Fix}\,\altmathcal{R}$; they are saddles, centers or, at bifurcations,  degenerate saddles. 
\item Outside of $\mathrm{Fix}\, \altmathcal{R}$, there can be only one pair of sink and source. They are always located in the synchrony subspace and therefore do not depend on the coupling strength $\kappa$. They are related  by $\altmathcal{R}$, which on the synchrony subspace induces also a time-reversal symmetry of the single uncoupled rotator.
\end{itemize}
The bifurcations of the equilibria will be discussed in detail below. 
The local bifurcations of the equilibria induce also changes of the configuration of the dissipative and conservative regions. Note that the regions can change also by global bifurcations given by a  reconnection of the saddle separatrices. The regions have the following general properties: \begin{itemize}
\item Each conservative region is filled with one-parametric family of neutral periodic orbits. 
\item Periodic orbits can have two different topological types: {\em Rotations}, where the curve closes after a full round trip of both oscillators such that both phases increase unboundedly, and {\em librations}, where both oscillators perform a small oscillatory motion whithout a full round trip in one of the phases. 
\item The conservative regions are bounded by homoclinics or heteroclinic cycles, which can also be of rotation or libration type. 
\item  Each dissipative region consists of heteroclinic orbits, connecting a source and a sink equilibrium.
\end{itemize}
The yellow and orange areas in Figs. \ref{Fig:Phase_Portr_1}--\ref{Fig:Phase_Portr_3}
indicate conservative regions filled with librations  (different colors correspond to clockwise and counter-clockwise motion), cyan and blue regions are filled with rotations.  
Dissipative regions can exist only in the presence of a source/sink pair of equilibria related by $\altmathcal{R}$ (white regions in  Figs. \ref{Fig:Phase_Portr_1}--\ref{Fig:Phase_Portr_3}). 
We give now a detailed description of the different types of local and global bifurcations occuring in this system. 

\paragraph*{Saddle-center bifurcation.} The red and magenta curves in the bifurcation diagram in Fig.~\ref{Fig: Bif_Diagr_S2} indicate a saddle-center bifurcation. At this bifurcation a saddle and a center equilibrium, both in $\mathrm{Fix}\altmathcal{R}$, merge and disappear, see also    Refs.~\cite{Broer_1998, Lerman2010,Lerman-Turaev2012}. The Jacobian at the degenerate equilibrium has an algebraically double zero eigenvalue. The resulting bifurcation condition for $a=1$ is given by 
\[
\omega=\pm\frac{\left(\sqrt{1+32\kappa^{2}}\pm 3\right)\sqrt{2\left(16\kappa^{2}
-1-\sqrt{1+32\kappa^{2}}\right)}}{32|\kappa|}
\]
and provides the red and magenta curves in Fig.~\ref{Fig: Bif_Diagr_S2}(a). An example of a phase portrait with such a degenerate equilibrium is given in Fig.~\ref{Fig:Phase_Portr_2}(k).  Together with the new equilibrium of center type there emerges also a conservative region, in this case filled with a family of periodic orbits of librations around this point. Also, a structurally stable homoclinic to the new saddle equilibrium emerges, giving the boundary of the conservative region, see e.g. the phase portrait in Fig.~\ref{Fig:Phase_Portr_1}(b). In Fig.~\ref{Fig:Phase_Portr_1} there are several pairs of structurally stable phase portraits related by this type of bifurcation: (b)--(a), (c)--(d), (e)--(f), (g)--h).
\paragraph*{Reversible pitchfork (sink/source) bifurcation.}
The green curves in Fig.~\ref{Fig: Bif_Diagr_S2}  indicate a reversible pitchfork bifurcation, where in a  spontaneous symmetry breaking a pair of a sink and a source equilibrium outside  $\mathrm{Fix}\altmathcal{R}$ bifurcate from an equilibrium within  $\mathrm{Fix}\altmathcal{R}$ that, at the same time, changes its type from a saddle to a center \cite{Lamb_1995, Lerman-Turaev2012}. The bifurcation condition is here a geometrically double zero eigenvalue. This bifurcation happens at the bifurcation of the single uncoupled oscillator at $\omega=\pm a$, where we have degenerate equilibria at  $(\phi_1,\phi_2)=(-\pi/2,\pi/2)$ and $(\pi/2,-\pi/2)$, respectively, that are  independent on $\kappa$.
Corresponding degenerate phase portraits are given in Figs.~\ref{Fig:Phase_Portr_2}(n),(o). The bifurcation of the equilibria induces also a reorganization of the homoclinic and heteroclinc connections and the conservative and dissipative regions. With the transformation of the center equilibrium in  $\mathrm{Fix}\altmathcal{R}$ into a saddle, the corresponding conservative region with periodic orbits surrounding the center vanishes. As  for the usual pitchfork bifurcation, the bifurcating pair of a source and sink equilibrium emerges with heteroclinic connections to the primary saddle equilibrium in  $\mathrm{Fix}\altmathcal{R}$. Moreover, as a consequence of the phase space being here the compact manifold $\mathbb{T}^{2}$, the  source and sink equilibrium inherit  global  heteroclinic connections to another saddle in  $\mathrm{Fix}\altmathcal{R}$, which, before the bifurcation, was carrying the homoclinic loop defining the boundary of the vanishing conservative region. Note that, obstructed by the new heteroclinic connections, also a conservative region of rotations vanishes. The pairs of structurally stable phase portraits related by this type of bifurcation in Fig.~\ref{Fig:Phase_Portr_1} are  (c)--(e) and (d)--(f). 
\paragraph*{Heteroclinic saddle-saddle connections.} 
We have two instances of  structurally unstable heteroclinic connections between  saddle equilibria in  $\mathrm{Fix}\altmathcal{R}$, given by the orange and black curves in the bifurcation diagrams in Fig.~\ref{Fig: Bif_Diagr_S2}.
The orange curves  indicate such a global bifurcation shown by the degenerate phase portraits in Figs.~\ref{Fig:Phase_Portr_2}(l),(m). This bifurcation induces the appearance/disappearance of conservative regions with rotations in between dissipative regions. Pairs of structurally stable phase portraits related by this type of bifurcation in Fig.~\ref{Fig:Phase_Portr_1} are  (b)--(c) and (a)--(d). 

\begin{figure}
\centering{}\includegraphics[width=0.5\textwidth]{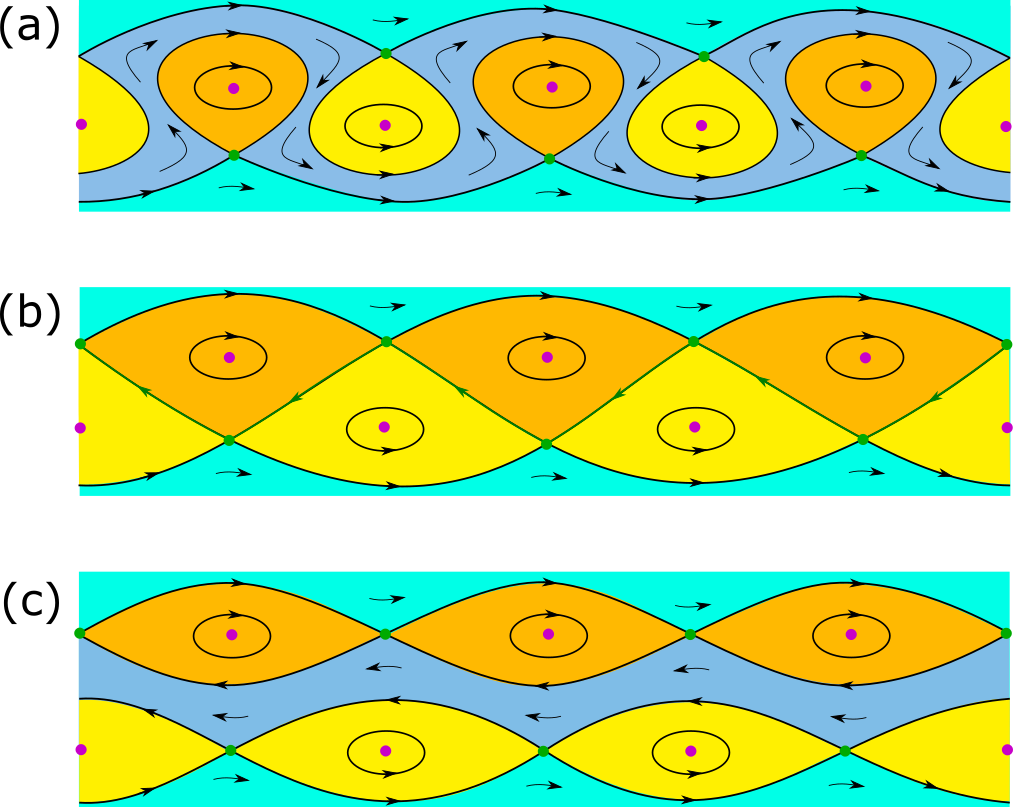}\caption{\label{Fig: Bif_Transit} 
Schematic phase portraits of the heteroclinic saddle-saddle connection (black curve in Fig.~\ref{Fig: Bif_Diagr_S2}(b)). 
Panel (a): structurally stable situation before the bifurcation (topologically equivalent to Fig.~\ref{Fig:Phase_Portr_1}(e) and Fig.~\ref{Fig: Phase_Portr_C2}(c)). 
Panel (b): degenerate situation at the bifurcation with structurally unstable saddle-saddle connection (topologically equivalent to Fig.~\ref{Fig:Phase_Portr_2}(p) and Fig.~\ref{Fig: Phase_Portr_C2cod-1}(j)). 
Panel (c): structurally stable situation after the bifurcation (topologically equivalent to 
 Fig.~\ref{Fig:Phase_Portr_1}(h) and  Fig.~\ref{Fig: Phase_Portr_C2}(d).
 Colors of invariant  regions and fixed points as in Fig.~\ref{Fig:Phase_Portr_1}.}
\end{figure}

The  black curves in Fig.~\ref{Fig: Bif_Diagr_S2}(b) correspond to heteroclinic saddle-saddle connections with a degenerate phase portrait as shown in Fig.~\ref{Fig:Phase_Portr_2}(p) and  connects the two structurally stable phase portraits (e) and (h)  in Fig.~\ref{Fig:Phase_Portr_1}. The mechanism how this global bifurcation leads to restructuring of the invariant regions is schematically shown in Fig.~\ref{Fig: Bif_Transit}. A region of undulating rotations (panel (a)) disappears and a new region of straight rotations in opposite direction appears (panel (c)). In the degenerate situation in between we see how two structurally stable homoclinics, which delineate the region of the rotations form two libration regions with opposite direction of motion, are reconnected through two  heteroclinic saddle-saddle connections forming a heteroclinic cycle of rotational type.

The bifurcation curves of such global bifurcation curves can typically be found only numerically. In our case of a planar flow this can be done by a simple shooting method. For the numerical treatment of more  general cases, see \cite{Champneys_1993}.
\paragraph*{Second time-reversal symmetry.} 
As already mentioned above, for $\omega=0$ the parametric symmetry $\gamma_1$ turns into a second time reversal symmetry $\altmathcal{R}_2$ with fixed space 
\begin{equation}
\mathrm{Fix}\altmathcal{R}_2=\left\{ (\phi_{1},\,\phi_{2}):\,\phi_{1}=\phi_{2}+\pi\right\}.
\label{eq: FixR2}
\end{equation}
This enables a homoclinic orbit to a saddle equilibrium  in  $\mathrm{Fix}\altmathcal{R}$ to turn into a heteroclinic connection between two saddle equilibria in  $\mathrm{Fix}\altmathcal{R}$ as soon as both saddles are related by $\altmathcal{R}_2$ and happens along the blue line in Fig.~\ref{Fig: Bif_Diagr_S2}(a). There are two qualitatively different degenerate phase portraits of this type given in Figs.~\ref{Fig:Phase_Portr_2}(i),(j), corresponding to the case of a  bifurcating  homoclinc orbit of libration and rotation type, respectively. Note that in the case of the rotation, there appears also a homoclinic to a saddle in $\mathrm{Fix}\altmathcal{R}_2$, which is structurally unstable with respect to perturbations that break the reversibility $\altmathcal{R}_2$. This type of global bifurcation mediates the transition of the structurally stable phase portraits (b) and (c) in Fig.~\ref{Fig:Phase_Portr_1} to their respective images under  $\altmathcal{R}_2$. 

For $\omega=0$ we find also two codimension-two bifurcations. For $a=\pm2 \kappa$ one of the two  equilibria in
$$\mathrm{Fix}\altmathcal{R}\,\cap\,\mathrm{Fix}\altmathcal{R}_2= \{(\pi/2,-\pi/2),(-\pi/2,\pi/2)\}$$
has a fully degenerate Jacobian, see
 phase portrait in Figs.~\ref{Fig:Phase_Portr_3}(q). At this point, two curves of saddle-center bifurcations that exist for $\omega\neq 0$ meet in a cusp point. Also the blue curve, indicating the heteroclinic connection between two saddle equilibria in  $\mathrm{Fix}\altmathcal{R}$, ends at this codimension-two point since the two saddles merge with the center in between them and vanish together with the enclosed conservative region.  
 
\begin{figure}
\centering{}\includegraphics[width=0.5\textwidth]{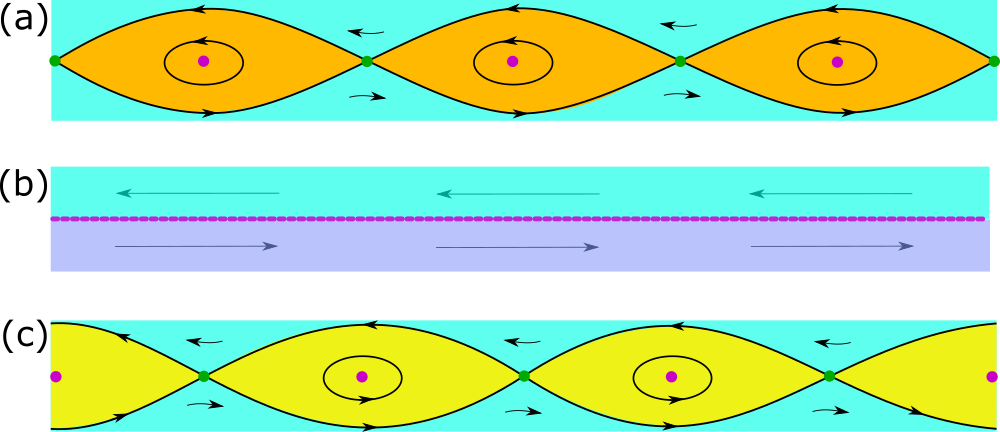}\caption{\label{Fig: Bif_Transit_Degen} 
Schematic phase portraits illustrating how locally the degenerate situation with rotational symmetry, 
occuring for $a=0$ (brown line in  Fig.~\ref{Fig: Bif_Diagr_S2}(b)), is unfolded for positive (panel(a)) and negative (panel (b)) values of $a$. 
The purple line in (b) consists of degenerate neutral fixed points.  Colors of invariant  regions and fixed points as in Fig.~\ref{Fig:Phase_Portr_1}. 
}
\end{figure}

 The second codimension-two bifurcation, with a degenerate phase portrait shown in  Figs.~\ref{Fig:Phase_Portr_3}(r), is a pair of saddle-saddle connections between a saddle in $\mathrm{Fix}\altmathcal{R}\,\cap\,\mathrm{Fix}\altmathcal{R}_2$ and a pair of saddles in  $\mathrm{Fix}\altmathcal{R}$ which are related by $\altmathcal{R}_2$. At this point meet two curves of saddle-saddle connections, which exist for $\omega\neq 0$.
%
\paragraph*{Rotational symmetry.} In the case of $a=0$ we have two coupled Kuramoto oscillators with phase shift symmetry, which can be reduced to a single equation for the phase difference. However, due to the special form of the coupling in the case (I) time-reversible system (\ref{eq:CaseI_1})--(\ref{eq:CaseI_2}), the phase difference $\psi=\phi_{1}-\phi_{2}$ stays always constant such that all trajectories are straight diagonal lines $\phi_{1}(t)=\phi_{2}(t)+\psi(0)$  with constant velocity $\dot \phi_{i}=\omega+\kappa\sin(\psi(0))$.
For $|\kappa|>|\omega|$ there are two diagonal lines $\phi_{2}=\phi_{1}\pm\arcsin(\omega/\kappa)$ 
with velocity zero. These lines of equilibria, existing along the brown bifurcation line in Fig.~\ref{Fig: Bif_Diagr_S2}(b), give rise to a global bifurcation in the following way. 
Breaking the phase shift symmetry with small $a\neq 0$, from each of the two lines of equilibria remains only a saddle and a center equilibrium while two narrow conservative regions of librations emerge, as shown in Fig.~\ref{Fig: Bif_Transit_Degen}. All other trajectories still form two regions of rotations with opposite directions, which are no straight lines any more but slightly modulated. 

The brown line ends in a codimension-two situation where the two lines of equilibria disappear together with region of rotations in opposite direction. The corresponding phase portrait is shown in Fig.~\ref{Fig:Phase_Portr_3}(t). Note that at this codimension two bifurcation points there emerge also curves of saddle-center bifurcations and heteroclinic saddle-saddle connections, see Fig.~\ref{Fig: Bif_Diagr_S2}(b). 
This bifurcation is described in more detail in Ref.~\cite{BMB_2022}.
\paragraph*{The uncoupled case.} For $\kappa=0$ the two rotators are decoupled. In this case the only bifurcation happens at $a=\pm \omega$, where both rotators simultaneously undergo the SNIC from excitable to rotating behavior. Note that the unfolding of this codimension-two point (see degenerate phase portrait in Fig.~\ref{Fig:Phase_Portr_3}(s)  for  $\kappa\neq 0$ gives rise to the reversible pitchfork bifurcation, where the stable synchronous equilibrium emerges, and, additionally, to curves of heteroclinic saddle-saddle connections and saddle-center bifurcations.
\paragraph*{Summary of case (I).}
Having clarified all the details of the bifurcation scenario in the time-reversible case (I) of two  rotators with anti-reciprocal coupling (\ref{eq:CaseI_1})--(\ref{eq:CaseI_2}) we can interpret the bifurcation scenario as follows. The transition of the single rotator from excitable to oscillatory motion at $|\kappa/\omega|=1$ plays the main role. If the single rotator is in the oscillatory regime, the coupled system has only conservative behavior. For weak coupling it consists of unidirectional rotation of both units. Stronger coupling leads to the coexistence of  bidirectional rotation and also to regions of libration, which can be seen as a conservative version of a {\em partial oscillation death}, i.e. for certain initial conditions the coupling prevents the rotating units from rotation or even reverses their rotation. 

For coupled rotators in the excitable regime we have always a pair of  source/sink equilibria that induces a dissipative region, which is identical to the basin of the sink. But only for small coupling this basin covers a set of full measure in phase space. For strong coupling there appear step by step conservative regions with both rotation and libration. In contrast to the oscillation death in the oscillatory regime, this can be seen as a {\em partial oscillation birth}, where -- again only for a  certain open set of initial conditions -- non-oscillatory units start to oscillate as a consequence of the coupling. 
%
%
%
%
\begin{figure}
\centering{}\includegraphics[width=\textwidth]{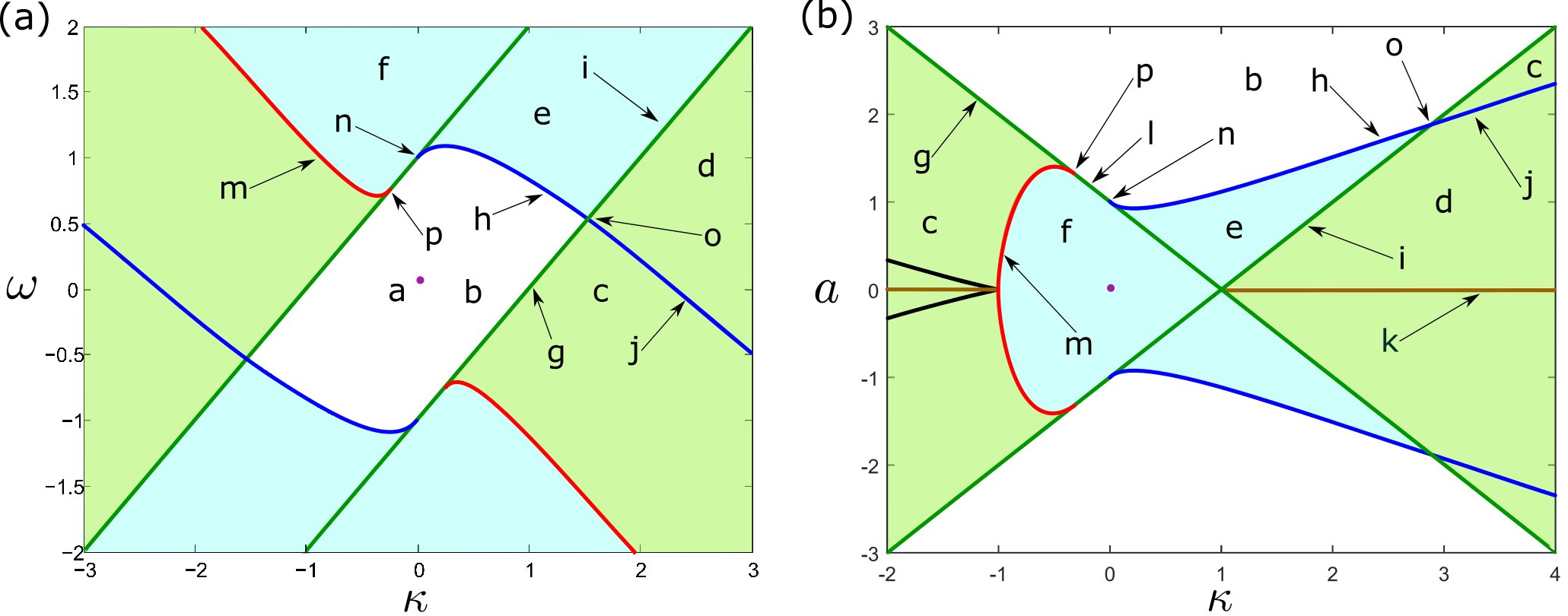}
\caption{
\label{Fig: Bif_Diagr_C2} 
Bifurcation diagrams for case (II) time-reversible system (\ref{eq:CaseII_1})--(\ref{eq:CaseII_2}).
(a) -- parameter plane $(\kappa,\omega)$  with fixed $a=1$. 
(b) -- parameter plane  $(\kappa,a)$ with fixed $\omega=1$.
Bifurcation curves: red -- saddle-center  bifurcation, 
green -- reversible equivariant sink/source bifurcation ,
blue/black -- heteroclinic saddle-saddle connections.
Moreover, there is a global bifurcation induced by the phase shift invariance (brown). 
Structurally stable phase portraits from the
different parameter regions with labels (a)–(f) are given in Fig.~\ref{Fig: Phase_Portr_C2},
phase portraits at the bifurcation curves and codimension-two points
(labels (g)–(p)) in Fig.~\ref{Fig: Phase_Portr_C2cod-1}. Colors of regions as in Fig.~\ref{Fig: Bif_Diagr_S2}.
}
\end{figure}
\subsection{Case (II): coupled rotators with reciprocal coupling}
\begin{figure}
\centering{}
\includegraphics[width=\textwidth]{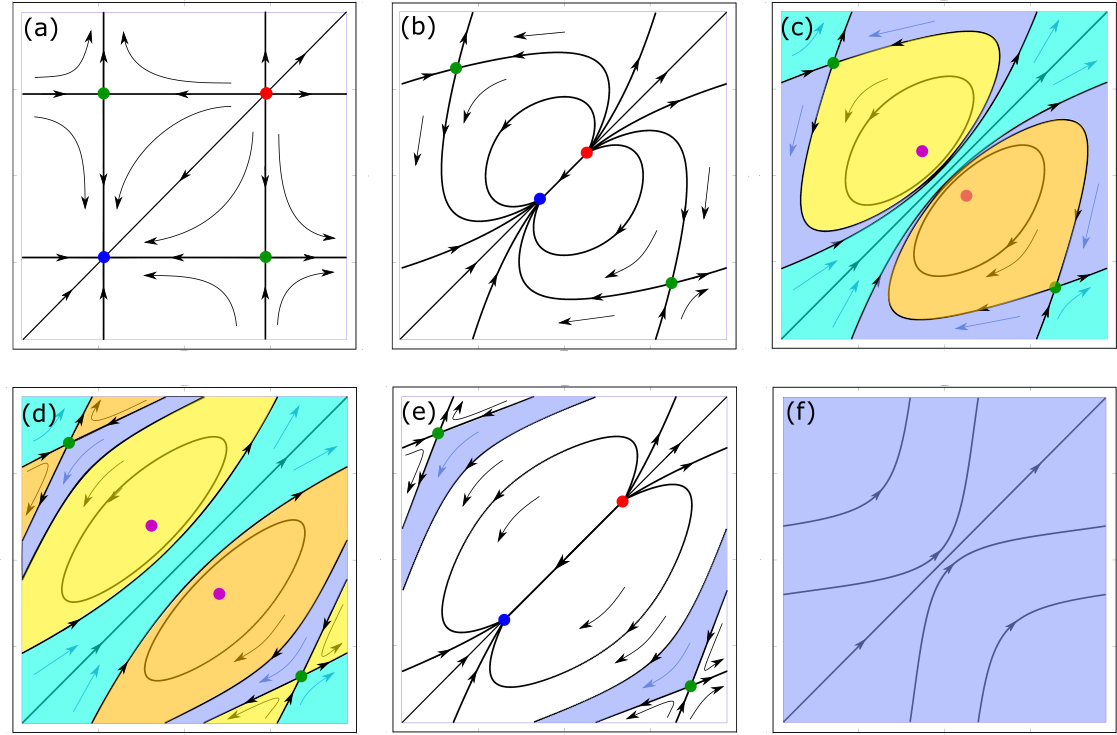}
\caption{\label{Fig: Phase_Portr_C2} Different types of structurally stable phase portraits in case (II) of time reversibility (\ref{eq:CaseII_1})--(\ref{eq:CaseII_2})  on the torus  $(\phi_{1},\phi_{2})\in[0,2\pi)\times[0,2\pi)$.
Parameters in panels (a)--(f)  are chosen from the correspondingly marked parameter regions in the bifurcation diagram Fig.~\ref{Fig: Bif_Diagr_C2}. 
Colors of invariant  regions and fixed points as Fig.~\ref{Fig:Phase_Portr_1}. 
}
\end{figure}

\begin{figure}
\centering{}
\includegraphics[width=\textwidth]{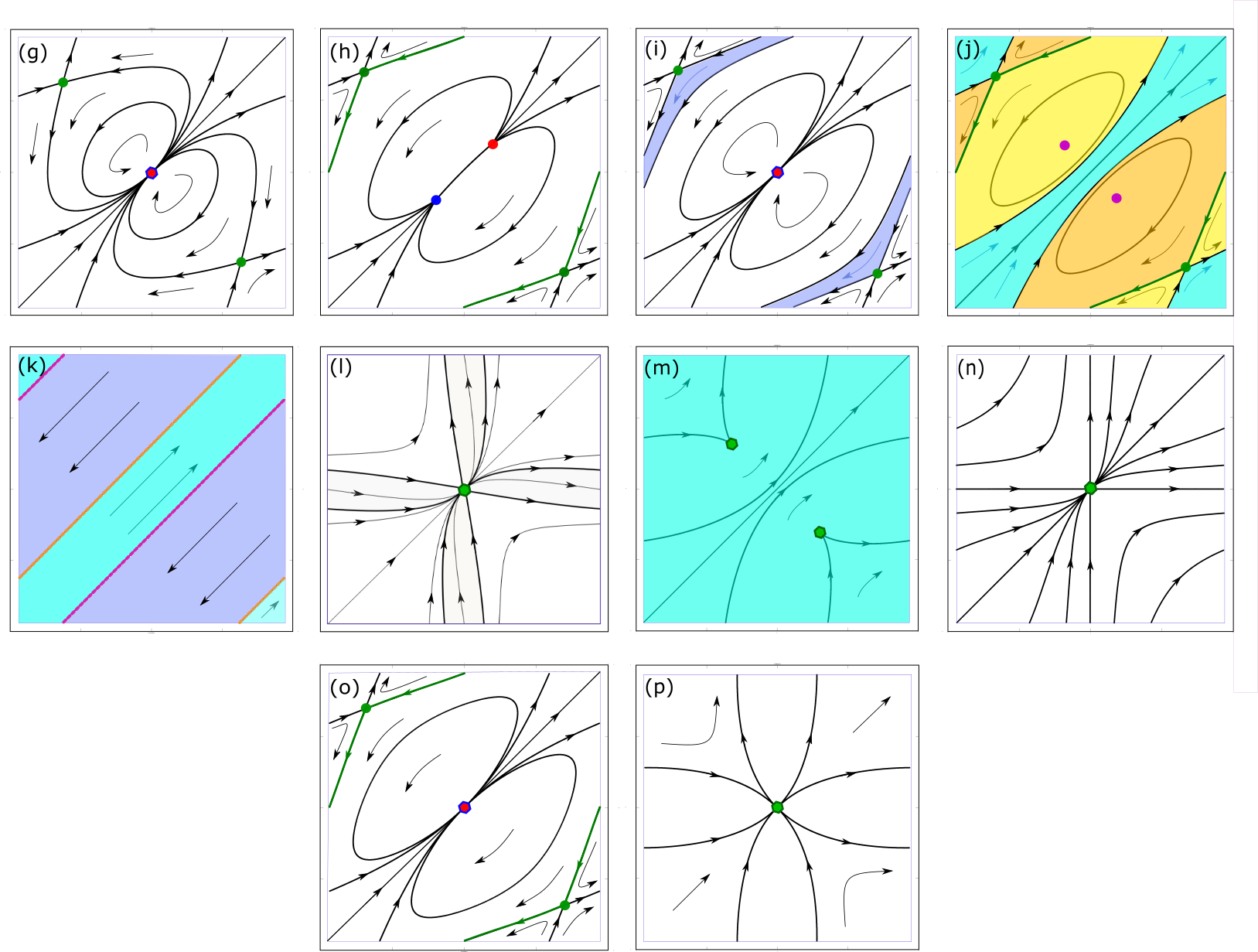}
\caption{\label{Fig: Phase_Portr_C2cod-1} Examples of structurally unstable phase portraits in case (II) of time reversibility (\ref{eq:CaseII_1})--(\ref{eq:CaseII_2})  on the torus  $(\phi_{1},\phi_{2})\in[0,2\pi)\times[0,2\pi)$. Panels (g)--(m): codimension-one; (n)--(p): codimension-two.  
Parameters are chosen from the correspondingly marked parameter regions in the bifurcation diagram Fig.~\ref{Fig: Bif_Diagr_C2}. 
Colors of invariant  regions and fixed points as Fig.~\ref{Fig:Phase_Portr_1}.
}
\end{figure}

 
We consider now the case (II) of time reversibility given by system (\ref{eq:CaseII_1})--(\ref{eq:CaseII_2}).
As above, we show in Fig.~\ref{Fig: Bif_Diagr_C2}  ~two bifurcation diagrams with respect to $(\kappa,\omega)$ with fixed $a=1$ (panel (a)) and with respect to $(\kappa,a)$ with fixed $\omega=1$  (panel (b)).
In addition to the time-reversal symmetry $\altmathcal{R}$, given by  (\ref{eq: R_Sym}), the system (\ref{eq:CaseII_1})--(\ref{eq:CaseII_2}) possesses the following $\mathbb{Z}_2$-equivariance
\[
\gamma_{m}:\ (\phi_{1},\,\phi_{2})\longmapsto(\phi_{2},\,\phi_{1}),
\]
which is a mirror symmetry with the invariant subspace 
\[
\mathrm{Fix}\gamma_{m}=\left\{ (\phi_{1},\,\phi_{2}):\,\phi_{1}=\phi_{2}\right\},
\]
corresponding to complete synchronization. Note that the composition
$$ \altmathcal{R}\gamma_{m}=\gamma_{m}\altmathcal{R}:\ (\phi_{1},\,\phi_{2},t)\longmapsto(-\phi_{1},\,-\phi_{2},-t)$$
of the time-reversal symmetry and the $\mathbb{Z}_2$-equivariance provides another time-reversal symmetry. This map has two invariant points $\{(0,0),(\pi,\pi)\}$, which will become important for the sink/source bifurcation discussed below.
Moreover, system (\ref{eq:CaseII_1})--(\ref{eq:CaseII_2}) has a time-reversal symmetry involving the parameters $\omega$ and $\kappa$:
\[
\gamma_{4}:\ (\phi_{1},\,\phi_{2},\,\kappa,\,\omega,\,\,t)\longmapsto(\phi_{1}+\pi,\,\phi_{2}+\pi,\,-\kappa,\,-\omega,\,\,-t).
\]
According to this symmetry the bifurcation diagram in Fig.~\ref{Fig: Bif_Diagr_C2}(a) is invariant under the reflection of both $\kappa$ and $\omega$ together.
The parametric symmetry
\[
\gamma_{3}:\ (\phi_{1},\,\phi_{2},a,\,\,t)\longmapsto(\phi_{1}+\pi,\,\phi_{2}+\pi,\,-a,\,t)
\]
that was already present in case (I) induces the reflection symmetry of the bifurcation diagram in Fig.~\ref{Fig: Bif_Diagr_C2}(b) with respect to $a$, while there is no reflection with respect to $\kappa$ alone, as in case (I). The system can have up to four fixed points, either two pairs of a saddle and a center, all in $\mathrm{Fix}\altmathcal{R}$, or two saddles in $\mathrm{Fix}\altmathcal{R}$ and a sink/source pair in the synchronization subspace $\mathrm{Fix}\gamma_{m}$.
We encounter here the similar types of reversible bifurcations as in the case (I) scenario described above. However, some of them are modified by the additional $\mathbb{Z}_2$-equivariance $\gamma_{m}$.
\paragraph*{Saddle-center bifurcation.}
The red line of the saddle-center bifurcation in Fig.~\ref{Fig: Bif_Diagr_C2}(a) is given by 
\[
\omega=-\frac{1}{8\kappa}-\kappa,\quad |\kappa|>\frac{1}{4}.
\]
According to the mirror symmetry, we have now a pair of symmetry related degenerate equilibria  $(\phi_{1}^{*},-\phi_{1}^{*})$ and $(-\phi_{1}^{*},\phi_{1}^{*})$, as seen in the degenerate phase portrait Fig.~\ref{Fig: Phase_Portr_C2cod-1}(m). Recall that together with the centers there appear also conservative regions of libration type, see Fig.~\ref{Fig: Phase_Portr_C2}(c). However, in contrast to case (I) where this bifurcation may induce the coexistence of conservative and dissipative regions, we find it here only in the purely conservative region, where it induces regions  librations in a fully rotating scenario Fig.~\ref{Fig: Phase_Portr_C2}(f).
\paragraph*{Reversible equivariant sink/source bifurcation.}
Similar as the reversible pitchfork bifurcation discussed in case (I), this bifurcation gives rise to a pair of sink/source equilibria outside $\mathrm{Fix}\altmathcal{R}$. However, in case (II) it
includes an interplay of the time reversal symmetry $\altmathcal{R}$ with the $\mathbb{Z}_2$-equivariance $\gamma_{m}$ and is  given by a degenerate equilibrium with  double zero eigenvalue that lies in 
$$\mathrm{Fix}\altmathcal{R}\cap \mathrm{Fix}\gamma_{m}=\{(0,0),(\pi,\pi)\}.$$
The corresponding bifurcation condition $\omega\pm a=\kappa$ provides the green lines in  Fig.~\ref{Fig: Bif_Diagr_C2}. This bifurcation comes here in two different versions that can not be distinguished on the linear level. The first type has a degenerate  phase portrait as given in Fig.~\ref{Fig: Phase_Portr_C2cod-1}(g). The degenerate equilibrium  connects a folded branch with two center equilibria in $\mathrm{Fix}\altmathcal{R}$, which are related by $\gamma_{m}$, with another folded branch in $\mathrm{Fix\,\gamma_{m}}$ containing a source and a sink equilibrium, which are related by $\altmathcal{R}$. The branches are organized as in a complex fold of the form $z^{2}+\mu=0$, $z\in\mathbb{C}$. The structurally stable phase portraits related by this type of bifurcation are Fig.~\ref{Fig: Phase_Portr_C2}(b) and (c). Note that this bifurcation connects a fully dissipative phase portrait with a fully conservative one. In the conservative phase portrait there are two further saddle equilibria in $\mathrm{Fix}\altmathcal{R}$, which each have gained in the conservative situation a structurally stable homoclinic orbit, delineating  two regions of librations around the center equilibria from two regions of rotation.

In the second type  has a degenerate  phase portrait as given in Fig.~\ref{Fig: Phase_Portr_C2cod-1}(l). The degenerate fixed point connects a pair of branches with saddle equilibria in $\mathrm{Fix}\altmathcal{R}$, which are related by $\gamma_{m}$, with the branch of the sink/source pair. In this situation both folded branches extend to the same side of the bifurcation such that all four involved equilibria coexist on one side of the bifurcation  and have all disappeared on the other side (Fig.~\ref{Fig: Phase_Portr_C2cod-1}(f)). 
The two types change at a codimension-two point on the green curve where the curve of saddle-center bifurcations (red) ends and the corresponding degenerate phase portrait is given in Fig.~\ref{Fig: Phase_Portr_C2cod-1}(p).
Note that there is a second codimension-two point along the green curve where another curve of  global bifurcations (heteroclinic saddle-saddle connections) ends.  This induces global change in the dissipative phase portraits emerging at the green line, showing beyond this point also a coexisting conservative region with rotations.

\paragraph*{Heteroclinic saddle-saddle connections}
We have two instances of structurally unstable heteroclinic saddle-saddle connections, given by the blue and black curves in the bifurcation diagrams in Fig.~\ref{Fig: Bif_Diagr_C2}.
At the blue curve a conservative region of backward rotations appears. Fig.~\ref{Fig: Phase_Portr_C2}(h) shows how this happens  from a purely  dissipative situation (Fig.~\ref{Fig: Phase_Portr_C2}(b)) where it leads to mixed-type dynamics (Fig.~\ref{Fig: Phase_Portr_C2}(e)). After the codimension-two point (Fig.~\ref{Fig: Phase_Portr_C2cod-1}(o)), the transition happens from a fully conservative situation (Fig.~\ref{Fig: Phase_Portr_C2}(c)), where the heteroclinic saddle-saddle connection (Fig.~\ref{Fig: Phase_Portr_C2cod-1}(j)) induces a second region of rotations in opposite direction (Fig.~\ref{Fig: Phase_Portr_C2}(d). 
This type of transition occurs also along the black curve and has been described schematically in  Fig.~\ref{Fig: Bif_Transit}. Note that 
this bifurcation occurs for $a=1$ only at very large values of $\omega$ such that it is out of the range of Fig.~\ref{Fig: Bif_Diagr_C2}(a).

\paragraph*{Rotational symmetry.} As in case (I), for $a=0$ we obtain two Kuramoto oscillators with a rotational symmetry. The  two intervals of the straight magenta line $a=0$, $|\kappa|\ge1$ in the bifurcation diagram in Fig.~\ref{Fig: Bif_Diagr_C2}(b) correspond to the
situation shown schematically in Fig.~\ref{Fig: Bif_Transit_Degen}, where the forced breaking of the rotational symmetry induces a global bifurcation and small libration regions emerge from a line of equilibria (Fig.~\ref{Fig: Phase_Portr_C2cod-1}(k)).

\paragraph*{Summary of case (II).} Similar to case (I), the dynamics in the synchronization subspace plays a central role for the dynamics. The SNIC bifurcation in this subspace, which comes here in the full system as the reversible equivariant sink/source bifurcation, induces together with the sink/source pair the dissipative dynamics. However, this bifurcation depends  here also on the coupling strength, i.e. it does not coincide with the SNIC of the uncoupled unit. In the dissipative regime large values of $\omega$ can lead to rotations coexisting with the dissipative region. But in contrast to case (I) there are no librations coexisting with a dissipative region.In the fully conservative regime we have again situations with rotations, librations, and additional rotations in opposite directions. Comparing the dynamics with and without coupling, we find again both situations, where the coupling enables rotations of excitable units (oscillation birth) or prevents rotations of rotating units (oscillation death). While this happens in most cases only for a part of the phase space, we have here also a case, where for increasing coupling two non-oscillating but excitable units make transition from a fully dissipative regime without any oscillations to a fully rotating regime. 
%
%

\section{Generic perturbations of the reversible cases}\label{sec:pert}
\begin{figure}
\centering{}\includegraphics[width=\textwidth]{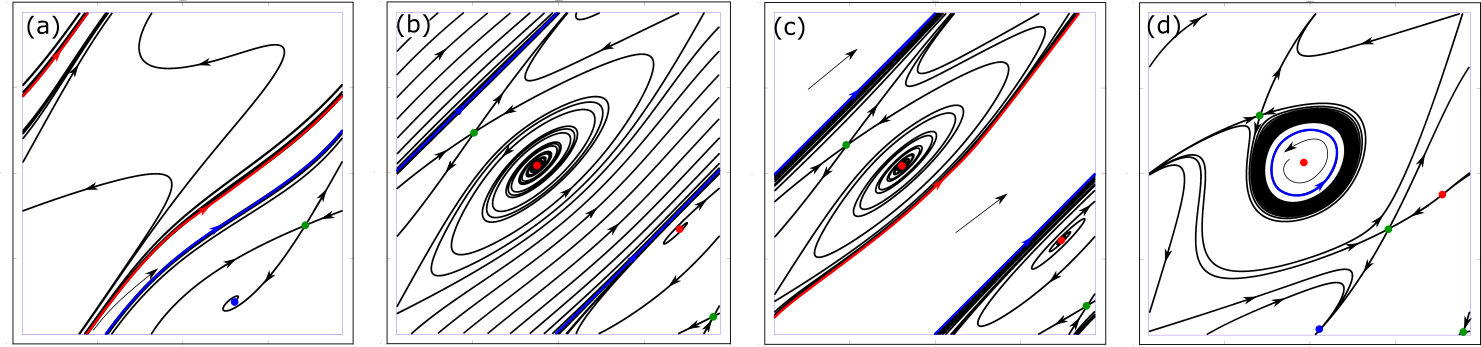}
\caption{\label{Fig: Ber_Asym}
Phase portraits of two coupled oscillators (\ref{eq:rot1})--(\ref{eq:rot2}) for different choices of the parameters $\omega_{1},\omega_{2},a_{1},a_{2},\kappa_{1},\kappa_{2}$:
(a) $(0.2,1,-1,0.5,3,2)$;
(b) $(0.5,0.5,0.21,0.2,1.1,1)$;
(c) $(0.7,0.7,0.1,0.1,1.06,1)$; 
(d) $(0.3,0.2,0.9,0.9,1.22,1.2)$.
Colors: red --- unstable, blue --- stable, green --- saddle.}
\end{figure}
For the general system (\ref{eq:rot1})--(\ref{eq:rot2}) of two coupled rotators, the reversible regimes studied above represent degenerate situations, which can be perturbed in different ways. Already for identical oscillators, i.e. $a_1=a_2, \omega_1=\omega_2$, and anti-reciprocal or reciprocal  coupling $\kappa_1=\pm\kappa_2$ a phase lag parameter $\alpha\neq 0$  or $\alpha\neq \pi/2$, corresponding to case (I) and case (II), respectively, will destroy the time-reversal symmetry. Other types of generic perturbations are non-identical oscillators or identical oscillators with different coupling strengths $\kappa_1\neq \pm\kappa_2$. 

Only the purely dissipative regimes close to the uncoupled non-oscillatory situation  have phase portraits that are structurally stable also under all such generic perturbations. As soon as there are conservative regions, all perturbations that break the time-reversal symmetry will lead to structural changes in the dynamics. 
\begin{itemize}
\item Center equilibria will turn into stable or unstable foci. 
\item There will be a slow drift along  families of neutrally stable periodic orbits in the conservative regions. 
\item Structurally stable homoclinic orbits, which constitute the boundaries of the conservative regions, will break. 
\end{itemize}
In this way, there can appear isolated stable and unstable periodic orbits from the families in the conservative regions. They can be both of rotation and libration type, see Fig.~\ref{Fig: Ber_Asym}. In particular, in the cases of perturbations of conservative dynamics with dissipative and conservative regions, this can lead to multistability, where new stable objects of different type emerge in addition to the structurally stable attracting equilibrium in the dissipative region as in panel (d). 
Note that for non-identical oscillators, there can be additional topologial types of  periodic orbits, where both units perform a different  number of round trips during one period. Such an example will be discussed in detail in the next section.

\subsection{Bursting-like orbits}
\begin{figure}
\centering{}\includegraphics[width=0.4\textwidth]{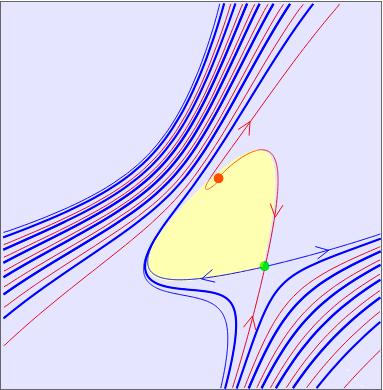}\caption{\label{Fig:burst1} Phase portrait for system (\ref{eq:rot1})--(\ref{eq:rot2}) with $a_{1,2}=1$, $\kappa_{1,2}=\mp1$ , $\alpha=0$ as in the reversible case (I), and $\omega_1=1.07, \,\omega_1=1.13$. Stable periodic orbit (blue), saddle (green) and unstable focus (red). Stable and unstable manifolds of saddle equilibrium (thin red and blue curves, respectively). Conservative regions of librations (yellow) and rotations (blue) from the reversible case with $\omega_1=\omega_2=1.1$}
\end{figure}
\begin{figure}
\centering{}\includegraphics[width=0.6\textwidth]{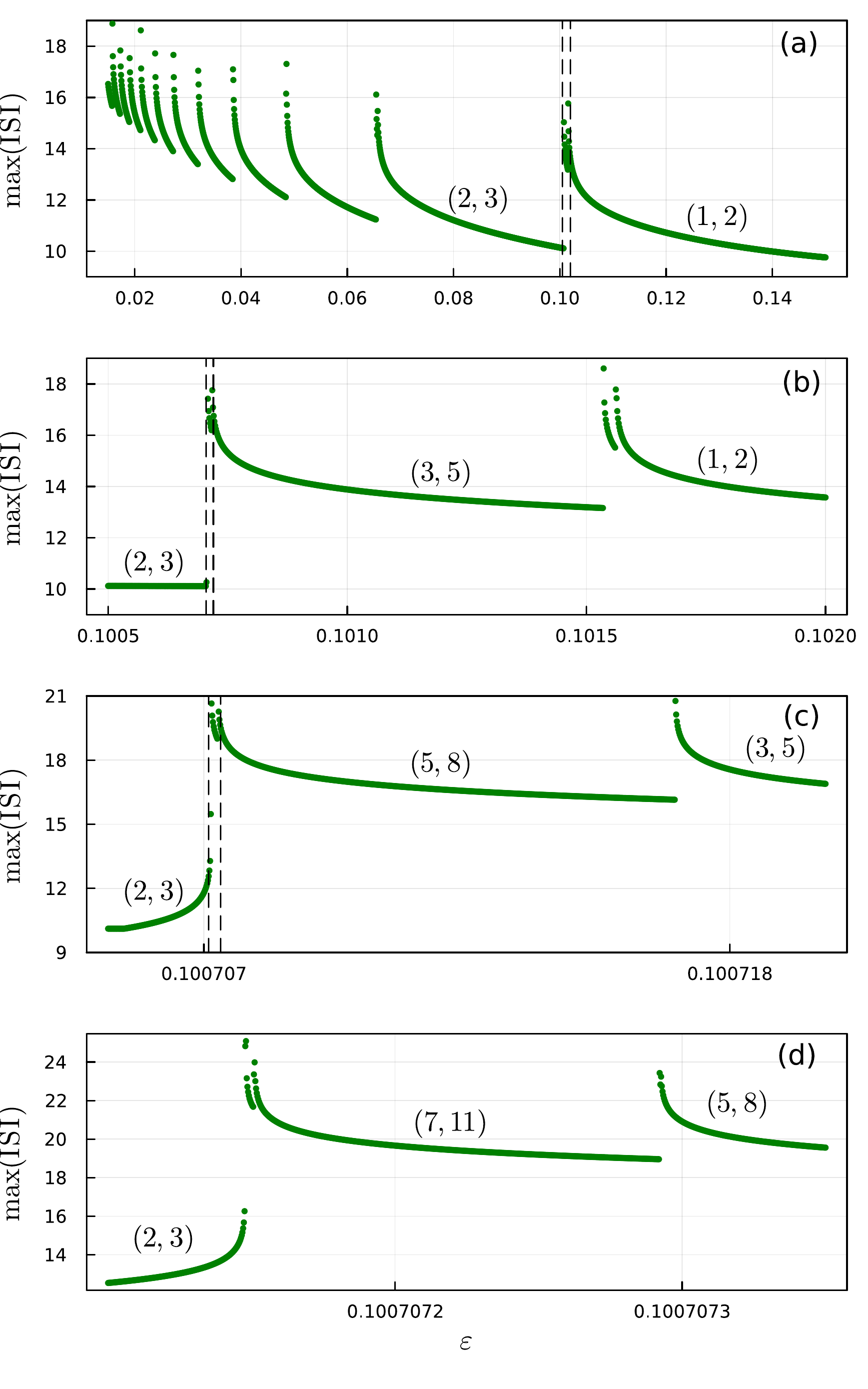}\caption{\label{Fig:ISI} Maximimum inter spike intervals for bursting solutions with varying detuning from the time reversible case for system (\ref{eq:rot1})--(\ref{eq:rot2}) with  $a_{1,2}=1$, $\kappa_{1,2}=\pm1$, $\omega_{1,2}=1.1\mp \varepsilon $, $\alpha=0$.}
\end{figure}
\begin{figure}
\centering\includegraphics[width=0.6\textwidth]{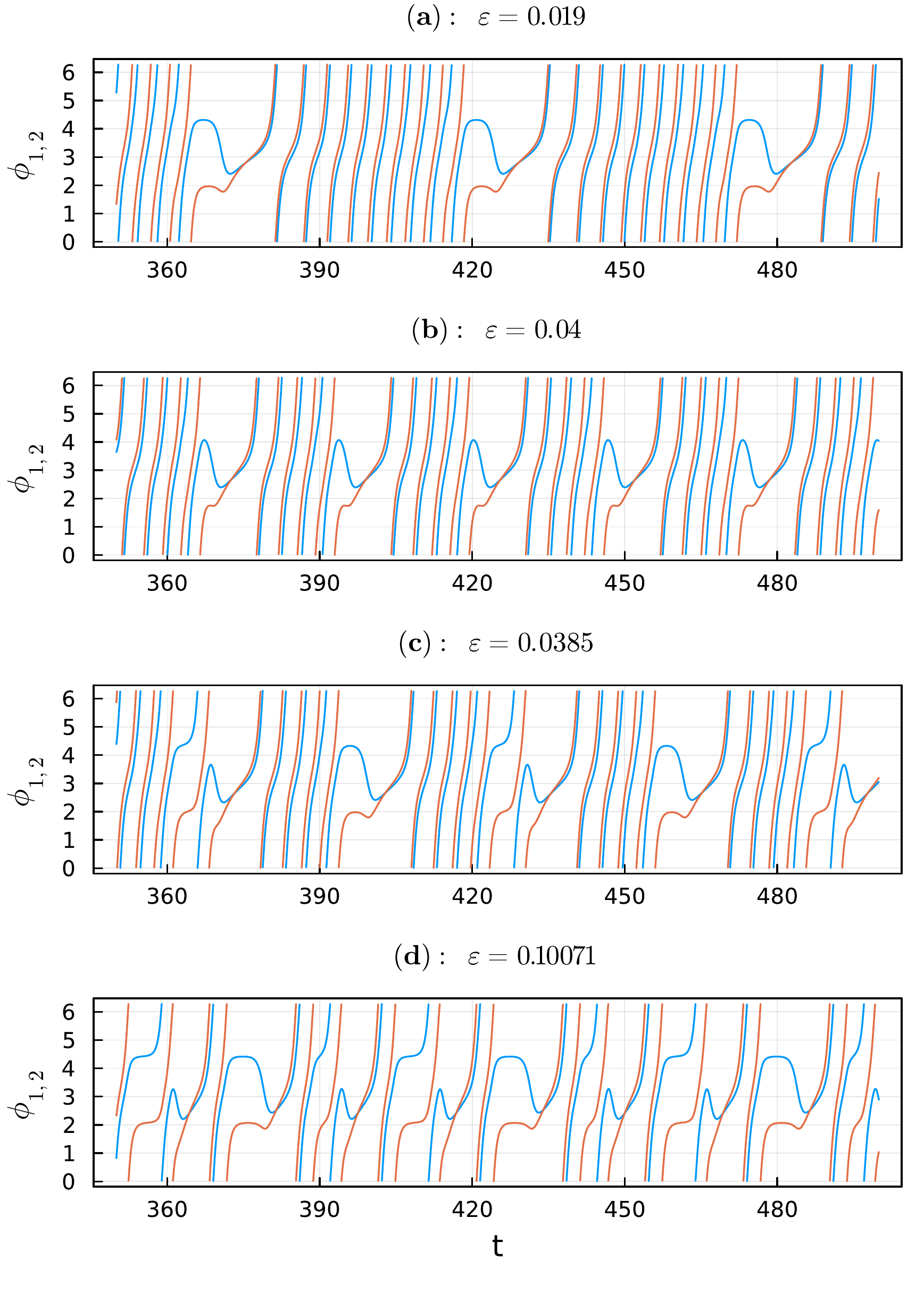}\caption{\label{Fig:burst3} Time traces of periodic solutions for system (\ref{eq:rot1})--(\ref{eq:rot2}); $\phi_1$--blue, $\phi_2$--red. For the indicated choices of the detuning value $\varepsilon$ we find the following winding numbers: (a)--$(9,10$), (b)--$(4,5)$, (c)--$(9,11)$, (d)--$(5,8)$. Other parameters as in Fig.~\ref{Fig:ISI}.
}
\end{figure}
A specific example of  non-trivial dynamics emerging from a small perturbation of the reversible dynamics in case (I), Fig.~\ref{Fig:Phase_Portr_1}(f) is shown in Figs.~\ref{Fig:burst1}--\ref{Fig:burst3}. In the time reversible case we have two conservative regions, one with librations and the other with rotations. 
A small perturbation to non-identical oscillators 
 with slightly detuned frequencies $\omega_1\neq\omega_2$
 induces a slow drift across the conservative region of rotations without stabilizing any of these rotations. At the same time, the center equilibrium within the other conservative region of librations is transformed into an unstable focus. The resulting  dynamics are shown in the  phase portrait in Fig.~\ref{Fig:burst1}.
We observe a stable periodic solution (blue) that performs a bursting-like behavior with many rotations during the slow passage through the conservative region of rotations until it comes close to the saddle equilibrium, where it can stay for an arbitrary long time interval.
For varying detuning of the frequencies the globally stable periodic orbits of this type  are organized in a complicated bifurcation scenario, where close to the conservative situation periodic solutions with arbitrarily long period and an increasing number of rotations within one burst appear.  In Fig.~\ref{Fig:ISI} we show how the branches of periodic solutions are organized  for varying the detuning $\varepsilon$ from the reversible  case at $\varepsilon=0$. Panel (a) shows a self-similar sequence of branches with increasing winding number $(n, n+1)$, $n=1,\dots,\infty$. Each of these branches ends at a homoclinic bifurcation, where the period grows unboundedly. However, at each of these transitions, we observe another self-similar cascade of transitions to orbits of more complicated structure existing only in increasingly small parameter windows.  Panel (b) shows the parameter region around the first transition in panel (a), where the branch with winding numbers $(1,2)$ disappears and a new branch with  $(2,3)$ appears. On the zoomed scale in panel (b) in between these two major branches a new branch with winding numbers $(3,5)$ becomes visible. Zooming into the transition between this branch and the $(2,3)$ branch, we find a branch with winding numbers $(5,8)$, see panel (c). Zooming in yet another time, we find a branch with winding numbers $(7,11)$, see panel (d). Examples of time traces of bursting orbits with different detuning values $\varepsilon$ and resulting winding numbers are given in Fig.~\ref{Fig:burst3}.
\subsection{Nonlinearities with higher harmonics}
\begin{figure}
\centering{}\includegraphics[width=\textwidth]{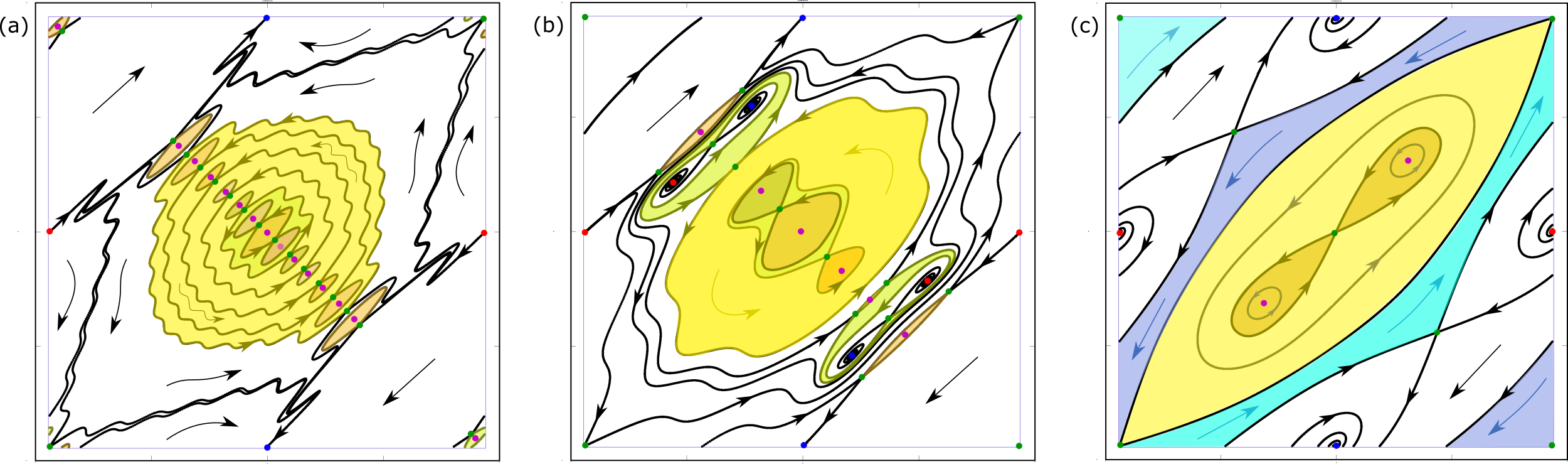}
\caption{\label{Fig: H_Harm} 
Examples of phase portraits for the general system of two
coupled active rotators (\ref{eq: General_1}), (\ref{eq: General_2}) with higher harmonics (\ref{eq:HH1}), (\ref{eq:HH2}). Parameters: 
(a): $\omega=0$, $p=0$, $\kappa=1.2$, $r=0.7$, $n=15$; 
(b): $\omega=0.1$, $p=0.2$, $m=3$, $\kappa=1.2$, $r=0.7$, $n=5$;
(c): $\omega=0$, $p=1$, $n=2$, $\kappa=4$, $r=0$.
Colors of regions and fixed points as in Fig.~\ref{Fig:Phase_Portr_1}.
}
\end{figure}
As explained in subsection \ref{subsec: Revers General}, the general  system (\ref{eq: General_1})--(\ref{eq: General_2}) has a time-reversal symmetry in two different cases, where
the functions $f_{1,2}(\phi)$, governing the local dynamics are identical and even, while the coupling functions $g_{1,2}(\phi)$ have to be also identical but can be either odd and have opposite signs (case (I)) or even and identical (case (II)). In section \ref{sec:dynamics-reversible} we investigated the case where both functions are restricted to the leading term in the Fourier expansion. In the general  case where both functions contain higher-order harmonics, the system can  possess more fixed points and, as a result, a much more complex structure of the invariant manifolds of the saddles, which provide the global structure of the dissipative and conservative regions in the regimes of mixed-type dynamics.

We will now briefly indicate how in the case (I) of a system with anti-reciprocal coupling the presence of higher-order harmonics in the local function of the local dynamics $f(\phi)=f_{1,2}(\phi)$ and in the coupling function $g(\phi)=g_{1}(\phi)=-g_{2}(\phi)$ can lead to more complex time-reversible dynamics.
First, note that already a single rotator $\dot\phi=f(\phi)$ with an  even function $f(\phi)$ containing the $n$-th harmonic $\cos(n\phi)$ can have up to $2n$ different fixed points. For a system of two such units with small coupling this gives rise to $4n^2$ fixed points, which are sinks, sources, and  saddles outside $\mathrm{Fix}\altmathcal{R}$ and also saddles inside $\mathrm{Fix}\altmathcal{R}$. At the other hand, for two  Kuramoto oscillators, i.e. $f(\phi)=\omega$, a coupling function $g(\phi)$ containing the $m$-th harmonic $\sin(\phi)$ can induce up to $2m$ lines of equilibria $\psi_j=\phi_1-\phi_2$, $j=1,\dots ,2m$, where $g(\psi_j)-\omega=0$. Breaking the rotational symmetry by a slightly non-constant $f(\phi)$, this leads to $2m$ saddle/center pairs in $\mathrm{Fix}\altmathcal{R}$ and corresponding conservative regions of rotations and librations, compare Fig.~\ref{Fig: Bif_Transit_Degen}. Hence, we can say that in a general system (\ref{eq: General_1})--(\ref{eq: General_2}) with  case(I) time reversibility, the  higher harmonics of 
$f(\phi)$ can induce multiple equilibria outside  $\mathrm{Fix}\altmathcal{R}$ and hence a more complex structure in the dissipative part, while functions  higher harmonics of $g(\phi)$ are responsible
for the emergence of multiple equilibria inside  $\mathrm{Fix}\altmathcal{R}$, leading to multiple conservative regions.

We illustrate this in Fig.~\ref{Fig: H_Harm} by two examples of functions of the form 
\begin{align}
    f(\phi)=&\omega-\cos\phi-p\cos(n\phi)\label{eq:HH1}\\
    g(\phi)=&\kappa(\sin\phi+r\sin(m\phi)).\label{eq:HH2}
\end{align}
Note that in panel (a) and (c) we have chosen $\omega=0$ such that we have the second time reversibility $\altmathcal{R}_2$, such that there are heteroclinic saddle-saddle connections, compare corresponding phase portraits in Figs.~\ref{Fig:Phase_Portr_2}(i),(j). In  panel (b), where we have chosen $\omega\neq 0$, all saddles in  $\mathrm{Fix}\altmathcal{R}$ have structurally stable homoclinics. Moreover, there are saddle equlibria outside  $\mathrm{Fix}\altmathcal{R}$. They come in pairs related by $\altmathcal{R}$ and can have structurally stable heteroclinic connections between them. They can be involved in a second type of reversible pitchfork bifurcation (saddle-saddle type), where an equilibrium inside $\mathrm{Fix}\altmathcal{R}$ changes from saddle to center, while a branch with two saddles outside $\mathrm{Fix}\altmathcal{R}$ emerges, cf Fig.~\ref{Fig: H_Harm}(b). 
Fig.~\ref{Fig: H_Harm}(c) shows  also  the result of a third type (center-center) of the reversible pitchfork bifurcations:  the emergence of two centers outside $\mathrm{Fix}\altmathcal{R}$ from a center inside $\mathrm{Fix}\altmathcal{R}$ that, at the same time, transforms from a center into a saddle. 

We see that some of the structural restrictions, which we encountered in the system (\ref{eq:CaseI_1})--(\ref{eq:CaseI_2}) of case (I) time-reversibility with only first harmonics are no more present for higher order harmonics. In particular,
\begin{itemize}
    \item there can be a large number of equilibria, in particular pairs of saddle equilibria outside $\mathrm{Fix}\altmathcal{R}$ and sink/source pairs outside the synchrony subspace.
    \item there can appear multiple nested regions of conservative regions of different type and nested regions of conservative and dissipative dynamics. 
\end{itemize}
However, the general observation remain true: The anti-reciprocal coupling of case (I) can induce a partial oscillation death, i.e. for certain initial conditions the coupling prevents the rotating units from rotation or even reverses their rotation. Also the effect of  partial oscillation birth, where non-oscillatory units start to
oscillate as a consequence of the coupling.

%
%

\section{Discussion and Outlook}
We have demonstrated that already for a fairly simple two-dimensional system of two coupled rotators in the transitional regimes between attractive and repulsive coupling there can arise quite complex dynamics. Particularly rich dynamics occur for parameter choices where the system has a time-reversal symmetry. In this case we also encounter the somewhat unusual types of  bifurcations of time-reversible systems. Additionally, such systems can  switch between  dissipative and  conservative dynamics, and also display the coexistence of different regions with such dynamics in phase space, which are governed by complex heteroclinic and homoclinic structures connecting the  fixed points within  and outside the symmetry subspace.

Note that certain interesting regimes and properties described in
this work for the system of two connected active rotators also exist
for more complex networks of rotators. 
In particular, a system of
$2N$ globally connected active rotators 
$\dot{\phi}_k=f_k(\phi_k)+\sum_{j=1}^{2N}g_{kj}(\phi_k-\phi_j)$,  where $f_k(x)=f_k(-x)=f_{k+N}(x)$, $g_{kj}(x)=\pm g_{kj}(-x)=\pm g_{k+N,j+N}(x)$,
can display a time-reversal symmetry similar to the cases given here, e.g. with a symmetry action of the form $\phi_{i}\mapsto-\phi_{i+N}$,
$i=1,\dots,N$. 
A system of $2N+1$ globally coupled active rotators with even coupling functions can have a time-reversal symmetry  with symmetry action $\phi_{i}\mapsto-\phi_{i+N+1}$, $i=1,\dots,N$, $\phi_{N+1}\mapsto-\phi_{N+1}$. 
Also  other symmetry actions, based on other permutations or phase shift symmetries, are possible. In all these cases, the system can have a coexistence of conservative and dissipative dynamics over wide regions in the parameter. Our preliminary numerical investigation indicates quite complex structures in a 4-dimensional system. 
We observe there conservative regions of multi-parameter families of neutral periodic orbits are bounded by  sets of homo/heteroclinic cycles. Despite the apparent complexity of global
bifurcations in multidimensional systems, certain of their properties
are similar to the bifurcations described above. Also the
destruction of conservative regions by small symmetry breaking perturbations and the emergence
of trajectories slowly drifting along the families of former neutral periodic orbits occurs in a somewhat similar way.


\section*{Acknowledgements}
Funded by the Deutsche Forschungsgemeinschaft (DFG, German Research Foundation) under Germany's Excellence Strategy – The Berlin Mathematics Research Center MATH+ (EXC-2046/1, project ID: 390685689).
O.B. acknowledges financial support of the Potsdam Institute for Climate Impact Research (PIK) and National Research Foundation of Ukraine (Project No. 2020.02/0089).
S.Y. was supported by the Deutsche Forschungsgemeinschaft (DFG, German Research Foundation), Project No. 411803875.


\vskip2pc

\bibliographystyle{RS}
\bibliography{Burylko,sy-references}

\end{document}